\documentclass[12pt,a4paper,centertags]{amsart}
\usepackage{pslatex}
\usepackage{hyperref}
\usepackage{graphicx}
\usepackage{pstricks}
\usepackage{mathrsfs}
\newtheorem{theorem}{Theorem}[section]
\newtheorem{lemma}[theorem]{Lemma}
\newtheorem{corollary}[theorem]{Corollary}
\newtheorem{proposition}[theorem]{Proposition}
\theoremstyle{remark}
\newtheorem{remark}[theorem]{Remark}
\theoremstyle{definition}

\numberwithin{equation}{section}

\DeclareMathOperator{\curl}{curl}

\newcommand{\loc}{{\text{\tiny loc}}}
\newcommand{\la}{\lambda}
\newcommand{\ga}{\gamma}
\newcommand{\al}{\alpha}
\newcommand{\e}{\mathtt{e}}
\newcommand{\im}{\mathtt{i}}
\newcommand{\ka}{\mathbf{k}}
\newcommand{\K}{\mathbf{K}}
\newcommand{\el}{\mathbf{l}}
\newcommand{\m}{\mathbf{m}}
\newcommand{\N}{\mathbf{N}}
\newcommand{\Z}{\mathbf{Z}}
\newcommand{\R}{\mathbf{R}}
\newcommand{\C}{\mathbf{C}}
\newcommand{\I}{\mathscr{I}}
\newcommand{\T}{\mathcal{T}}
\newcommand{\Pb}{\mathsf{P}}
\newcommand{\E}{\mathsf{E}\,}
\newcommand{\sF}{\mathscr{F}}
\newcommand{\sT}{\mathscr{T}}
\newcommand{\uno}[1]{\mathbf{1}_{#1}}
\newcommand{\pde}{\textsf{PDE}}
\newcommand{\ode}{\textsf{ODE}}
\newcommand{\app}{{\textsf{\Tiny (apprx)}}}
\newcommand{\esp}{{\textsf{\Tiny (expl)}}}
\newcommand{\pz}{{\scriptscriptstyle(0)}}
\newcommand{\pn}{{\scriptscriptstyle(n)}}
\newcommand{\pnm}{{\scriptscriptstyle(n-1)}}
\begin{document}
\title[A branching process representation...]{A probabilistic
representation for the solutions to some non-linear PDEs using 
pruned branching trees}
\author{D. Bl\"omker}
\address{Institut f\"ur Mathematik, RWTH Aachen, Templergraben 55, D-52052
Aachen, Germany}
\email{bloemker@instmath.rwth-aachen.de}
\author{M. Romito}
\address{Dipartimento di Matematica "U. Dini", viale Morgagni 67/a, I-50134
Firenze, Italia}
\email{romito@math.unifi.it}
\author{R. Tribe}
\address{University of Warwick, Mathematics Institute, CV4 7AL Coventry, UK}
\email{tribe@maths.warwick.ac.uk} 
\begin{abstract}
The solutions to a large class of semi-linear parabolic {\pde}s are given
in terms of expectations of suitable functionals of a tree of
branching particles. A sufficient, and in some cases necessary,
condition is given for the integrability of the stochastic representation, 
using a companion scalar {\pde}.

In cases where the representation fails to be integrable
a sequence of pruned trees is constructed, producing a
approximate stochastic representations that in some cases converge, 
globally in time, to the solution of the original {\pde}.
\end{abstract}
\maketitle
%%
%%%%%%%%%%%%%%%%%%%%%%%%%%%%%%%%%%%%%%%%%%%%%%%%%%%%%%%%%%%%%%%%%%%%%%%%%%%%%%%%%
%%
%%
\section{Introduction}
This paper considers stochastic representations for solutions to  
a large class of semi-linear parabolic {\pde}s, or systems of {\pde}s, of the type 
\begin{equation}\label{e:main}
\partial_t u = Au + \sF(u) + f,
\end{equation}
where $A$ is an operator with a complete 
set of eigenfunctions, $\sF$ is a polynomial nonlinearity
in $u$ and its derivatives, and $f$ is a given driving function.

In short, the solution $u$ is expanded into a Fourier series 
using the eigenfunctions of $A$.  This yields (as in 
spectral Galerkin methods) a system of 
countably many coupled {\ode}s for 
the Fourier coefficients. 
This {\ode} system is then solved in a weighted 
$\ell^\infty$-space, 
via an expectation over a tree of branching particles. 
The rules for the branching and dying probabilities 
arise from the particular {\pde} being studied. Moreover the {\pde} 
determines an evaluation operator $R_t$, acting on the tree $\T_\ka$ of 
particles rooted at each Fourier mode $\ka$, so that, under integrability
assumptions, the (suitably weighted) $\ka$'th Fourier mode $\chi_{\ka}(t)$
is given by 
$$ 
\chi_\ka(t)=\E[ R_t(\T_\ka)].
$$
This is precisely the method of Le Jan and Sznitman \cite{LJS}
(later extended in Bhattacharya et
al.\ \cite{BCDG}, Waymire \cite{Way}, or Ossiander \cite{Oss})
where they treated the Navier-Stokes equations in $\R^3$.
Earlier papers connecting branching particle systems to {\pde}s 
(for instance Skorokhod \cite{Sko} 
or Ikeda, Nagasawa, and Watanabe \cite{INW},
and later McKean \cite{McK})
use branching coupled with a diffusion,
and the stochastic representation  is derived directly without 
Fourier series, so that the linear operator 
$A$ is limited to generators of diffusions.

Our first aim is to show that this method applies to a large range of 
equations. In Section \ref{abstract} we present three representative examples 
with quadratic non-linearities. We comment on further generalisations 
in Section \ref{ss:moregen}. Basically any system of parabolic {\pde}s 
with polynomial nonlinearities in the derivatives is admissible.
 
One major drawback of the stochastic representation is that  it often
fails to exist for large times $t$,  although the solution to the
{\pde} may still exist.  The problem is that $R_t(\T_\ka)$ may fail to be
an integrable random variable for $t \geq t_0$.   
When  deriving a system of {\ode}s in $\ell^{\infty}$ space, there
is considerable freedom  in the choice of  weights for the Fourier
coefficient under which integrability can be established.  See
Bhattacharya et al.\ \cite{BCDG} for an extensive discussion in the
case of  $3D$-Navier Stokes. 
One way to check
integrability is to establish a scalar comparison  equation. The
finiteness of the comparison equation implies the integrability needed
for the stochastic representation to hold.  
The comparison equation is
independent of the weights and represents a worst case scenario
with super-linear (explosive) growth. It typically completely
ignores most of the structure of the non-linearity in the
original {\pde}.

In Section \ref{s:rep} we establish the stochastic 
representation under the assumption that it is integrable. 
In Section  \ref{compare} we investigate the comparison equation.
This typically shows the representation is integrable at small times, or,
when there is no linear instability, 
for all times with small data. For some classes of equations, for 
example $1d$-Burgers equation, 
we obtain a necessary and sufficient condition for the 
integrability of the stochastic representation,
independent of the choice of weights in Fourier space.

Our second aim is to present an approach to treat cases
where integrability fails. The key point is to get rid of the smallness
condition on the initial data and the forcing, in order to
find a stochastic representation that is global in time. 
In Bhattacharya et al.\ \cite{BCDG} the branching trees are 
pruned after $n$ generations. This gives a stochastic representation 
of a Picard iteration scheme converging to the original 
{\pde}, but, as stated in \cite{BCDG}, 
the existence of the expectation is equivalent to the convergence of 
the  Picard iteration scheme. In another approach, Morandin \cite{Mor}
suggested a clever re-summation of the expectation in order to
improve the convergence for large times, but he was only able to
rigorously verify the global convergence of his method in a
simple example where \eqref{e:main} is a one-dimensional {\ode}.
Our approach is to construct sets $\Omega_n$, with $\Pb[\Omega_n] \uparrow 1$
so that
$$
\chi_{\ka}(t) = \lim_{n \to \infty} \E \left[ R_t(\T_\ka) \, 
\uno{\Omega_n} \right].
$$
This treats the expectation somewhat as a singular integral, where we
have to be careful how to cut out the singularity. 

The method we use, explained in Section \ref{sec4}, is to construct 
a pruned branching tree  $\T_\ka^\pn$
which will agree with $\T_\ka$ on $\Omega_n$. The expectation
for the pruned tree
$$
\E[ R_t(\T_\ka^\pn)] = \E[R_t(\T_\ka)\,\uno{\Omega_n}]
$$ 
is well defined and represents the $\ka^\text{th}$ Fourier mode of the
solution to a semi-implicit approximation scheme of the type 
$$
\partial_t u^\pn
=A u^\pn + \tilde{\sF}(u^\pn,u^\pnm) + f.
$$
We then use {\pde} techniques to verify that the 
approximation scheme converges to a solution of the original {\pde}.
Although there are general results for the convergence of such
approximations (cf. for example Bj{\o}rhus and Stuart \cite{BjSt}) 
the assumptions are usually quite restrictive.
Since stronger arguments are model specific, we present the 
arguments only in two special cases, namely for 
a simple quadratic {\ode} and for Burgers equation. 
We believe these examples illustrate that the method potentially
dramatically extends the range of {\pde}s for which there is a global
stochastic representation. 
A global result is essential if one wants to study a stochastic 
representation of the long time behaviour of solutions,
for example in terms of stationary solutions or pull-back fixed points.
Only in a very simple framework of small initial conditions and 
uniformly small forcing is it currently possible to derive such results
(for the  $3$D Navier-Stokes
equations see Bakhtin \cite{Bak} and Waymire \cite{Way}).
The extension of such results to non-trivial cases
and the relation with the pruned representation
are the subject of work in progress.
%%
%%
%%%%%%%%%%%%%%%%%%%%%%%%%%%%%%%%%%%%%%%%%%%%%%%%%%%%%%%%%%%%%%%%%%%%%%%%%%%%%%%%%
%%
%%
%%
\section{Abstract setting and examples} \label{abstract}
We first present an infinite system of {\ode}s involving a quadratic
non-linearity. The system is indexed over $\ka\in\Z^d$. 
We then discuss several examples of {\pde}s on the torus $[0,2 \pi)^d$
and recast their Fourier transforms into our abstract {\ode} setting.  
We do not present the highest generality possible, but focus instead only 
an equation with one quadratic nonlinearity, one linear instability and
one forcing term. We comment in Subsection \ref{ss:moregen} on a large 
number of possible extensions, including other domains and boundary
conditions, multiple forcing terms and additional nonlinearities,
possibly of higher order.
\subsection{The general system of {\ode}s}
We consider solutions $\chi(t):\Z^d \to \C^r$ to the following
infinite dimensional system of $\C^r$-valued {\ode}s
\begin{equation}\label{geneq}
\dot\chi_\ka = \lambda_\ka \Big[ -\chi_\ka
   + C_f p_\ka \chi_\ka
   + C_b \sum_{\el,\m\in\Z^d} q_{\ka,\el,\m} B_{\ka,\el,\m}(\chi_{\el},\chi_{\m})
   + d_\ka  \gamma_\ka \Big].
\end{equation}
with $\ka\in\Z^d$. The constants $\lambda_\ka>0$ (which will determine 
the rate
of particle evolution), $p_\ka$, $q_{\ka,\el,\m}$, $d_\ka\in[0,1]$
(which will determine the probabilities of flipping, branching and dying),
and $C_f$, $C_b\ge0$ (the flipping and branching constants) are
fixed, as are bilinear operators 
$B_{\ka,\el,\m}:\C^r \times \C^r \to\C^r$ satisfying  
$$
|B_{\ka,\el,\m}(\chi,\chi')| \leq |\chi| \, |\chi'|
$$
for all $\chi,\chi' \in \C^r$. The choice of these constants will
arise from the Fourier transform of the {\pde} being studied. We
assume throughout that
\begin{equation} \label{e:probsum}
p_\ka + q_{\ka} + d_\ka=1 \qquad \text{for all }\ka \in \Z^d,
\end{equation}
and 
\begin{equation} \label{asslocex}
p_\ka \to 0,\quad q_{\ka}\to 0, \quad \text{as }|\ka|\to\infty
\end{equation}
where
$$
q_{\ka} = \sum_{\el,\m \in Z^d} q_{\ka,\el,\m}.
$$
The data for the equations consists of a time dependent forcing 
$\ga = \{\ga_\ka(t): k \in \Z^d, \, t \geq 0 \}$ and an 
initial condition $\chi(0) = \{\chi_\ka(0): \ka \in \Z^d\}$. 
We consider the above system in its mild formulation, that is for
given data we look for measurable $t\mapsto\chi_\ka(t)\in\C^r$
satisfying, for $\ka \in \Z^d$,
\begin{multline}\label{mildform}
\chi_\ka (t)
= \e^{-\lambda_\ka t}\chi_\ka(0)
 +\int_0^t\lambda_\ka \e^{-\lambda_\ka(t-s)}
   \Big[ C_f p_\ka \chi_\ka(s)+\\
    +C_b \sum_{\el,\m\in\Z^d} q_{\ka,\el,\m} B_{\ka,\el,\m}(\chi_\el(s),\chi_\m(s))
    +d_\ka\gamma_\ka(s) \Big]\,ds.
\end{multline}
Note that we need some regularity of $\chi_\ka$, in order to make
\eqref{mildform} well defined.
\begin{remark} \label{Remark2.1}
There is considerable flexibility when choosing the 
constants in the {\ode} system \eqref{geneq}. 
For example, we can adjust the probabilities $p_\ka$, $q_{\ka,\el,\m}$, and 
$d_\ka$ by adjusting the constants $C_b$, $C_f$ and considering modified
forcing data $\ga$. In particular, in an equation where the probabilities
do not add up to $1$ in \eqref{e:probsum}, it is always possible 
to adjust $d_\ka$ and the forcing data so that this constraint holds. 
Similarly, an equation with $C_f$ and $C_b$ replaced by bounded functions
of $\ka$ can be recast into the form \eqref{geneq} by forcing the $\ka$ 
dependence into the probabilities $p_{\ka}$, $q_{\ka}$, and $d_{\ka}$.
\end{remark}
\subsection{The d-dimensional Burgers equations}\label{suse:burg}
Consider solutions $u(t,x) \in \R^d$, for $t \geq 0$ and $x\in[0,2\pi)^d$
to the Burgers system
\begin{equation}\label{burgers}
\begin{cases}
\partial_t u-\Delta u+(u\cdot\nabla)u=f,\\
u(0)=u^0,
\end{cases}
\end{equation}
with periodic boundary conditions, where $f$ is an external forcing. 
We restrict ourselves to periodic boundary conditions, 
as the nonlinearity is easy to compute in the Fourier basis. Nevertheless, 
other kinds of boundary conditions, for instance like Dirichlet or Neumann,
can be treated in a similar fashion (cf. section \ref{ss:moregen}).

If we expand the solution 
$$
u(t,x)=\sum_{\ka\in\Z^d}u_\ka(t)\e^{\im\ka\cdot x},
$$
the equation reads in the Fourier coefficients as
$$
\dot u_\ka=-|\ka|^2u_\ka-\im\sum_{\el+\m=\ka}(u_\el\cdot\m)u_\m+f_\ka,
$$
The sum is over all $\el, \m\in\Z^d$ satisfying $\el+\m=\ka$. 
Define a weight function $w_{\ka} = 1 \vee |\ka|^{\alpha}$, where
$\alpha>0$ will be chosen shortly, and set $\chi_\ka = w_\ka u_\ka$.
Then
\begin{equation}\label{e:burg-F}
\begin{cases}
\dot{\chi}_\ka =-|\ka|^2 \chi_\ka-\im \sum_{\el+\m=\ka}\frac{|\m|w_\ka}{w_\m w_\el}(\chi_\el\cdot
\frac{\m}{|\m|}) \chi_\m + f_\ka w_\ka,\\
\chi_\ka(0) = u_\ka(0) w_\ka\;.
\end{cases}
\end{equation}
Note that the mode $u_\mathbf{0}$ has no linear dissipation. Below 
we will add and subtract $\la_\mathbf{0}\chi_\mathbf{0}$ to the
equation for $\chi_\mathbf{0}$, which introduces a linear instability
but which allows us to write the equation in our desired abstract form.
We note that in dimension $d=1$ this trick is unnecessary:  
the equations for the zero${}^\text{th}$ mode decouples, in that it
simplifies to $\dot u_\mathbf{0} = f_\mathbf{0}$, and it is then possible
to reduce the problem to the case $f_\mathbf{0}= u_\mathbf{0} = 0$.

We now show one way to recast \eqref{e:burg-F} into the 
abstract form \eqref{geneq}. For given $C_f$, $C_b$, $\la_\mathbf{0}>0$
we define
\begin{gather*}
\lambda_\ka = \begin{cases}|\ka|^2 &\ka\neq\mathbf{0},\\
              \la_\mathbf{0} &\ka=\mathbf{0},\end{cases}
\qquad 
p_\ka = \begin{cases}0 &\ka\neq \mathbf{0},\\
        C_f^{-1} &\ka=\mathbf{0},\end{cases}\\
q_{\ka,\el,\m} = C_b^{-1} \frac{|\m|w_\ka}{\la_\ka w_\el w_\m}, 
\qquad 
B_{\ka,\el,\m}(\chi, \chi') = - \im (\chi \cdot \frac{\m}{|\m|}) \chi',
\end{gather*}
whenever $\el + \m = \ka$ (and zero otherwise).
Lemma \ref{alphamajo} below ensures, provided 
we choose $\alpha > \max \{\frac{d+1}{2},d-1\}$, that 
$q_\ka = \sum_{\el, \m} q_{\ka,\el,\m} < \infty $ and that
$q_\ka \to 0 $ as $|\ka| \to \infty$. Thus by taking $C_b$,
$C_f$ sufficiently large we have that $p_\ka+q_\ka < 1$
and it remains only to define $d_\ka = 1-p_\ka-q_\ka$ and 
$\ga_\ka = (f_\ka w_\ka/\la_\ka d_\ka) $ for $\ka \in \Z^d$.
Again,  there is considerable flexibility in these choices.
\begin{lemma}\label{alphamajo}
For all $\alpha$, $\gamma >0$ with $\alpha+\gamma>d$ there exists
$C = C(\alpha, \gamma) < \infty$ so that, for all $\ka \in \Z^d$,
with $\ka\neq\mathbf{0}$,
$$
\sum_{\substack{\el+\m=\ka\\\el\neq\mathbf{0},\m\neq\mathbf{0}}}
\frac1{|\m|^\alpha|\el|^\gamma}
\le\begin{cases}
C\left( 1+ |\ka|\right) ^{-\beta},&\qquad\text{if $\alpha\neq d$ and $\gamma\neq d$},\\
C \left(1+|\ka|\right)^{-\beta} \log(1+|\ka|),
&\qquad\text{if $\alpha=d$ or $\gamma=d$},\\
\end{cases}
$$
where $\beta=\min\{\alpha, \gamma, \alpha+\gamma-d\}$ and the sum is over
all indices $\el$, $\m$ in $\Z^d$ satisfying the given constraints.
\end{lemma}
One way to prove this lemma, whose proof is omitted, 
is to compare above and below by suitable
continuous integrals. 
\subsection{Two dimensional Navier-Stokes equations}\label{navierstokes}
We briefly treat the two dimensional Navier-Stokes in its
vorticity formulation, since this will be used 
in section \ref{compare} as an example where the 
comparison equation yields exact statements about 
the integrability of the stochastic representation. 

In dimension $d=2$ the vorticity $\xi=\curl u$ is a scalar
and satisfies, on the torus $[0,2\pi)^2$ and with periodic 
boundary conditions,  
\begin{equation}\label{vorticity}
\begin{cases}
\partial_t\xi-\Delta\xi+(u\cdot\nabla)\xi=f,\\
\xi(0)=\xi^0,
\end{cases}
\end{equation}
where $u$ is the solution to the Navier-Stokes equations. 
The Fourier coefficients satisfy the following system, 
$$
\dot\xi_\ka = - |\ka|^2\xi_\ka + \sum_{\el+\m=\ka}
\frac{\ka\cdot\el^\perp}{|\el|^2}\xi_\el\xi_\m + f_\ka,
$$
where $\el^\perp=(l_2,-l_1)$. For simplicity we shall assume 
that $f_\mathbf{0} =0$ and the vorticity has mean value
$\xi_\mathbf{0}$ zero and is omitted from the system. 

We then set $\chi_\ka=|\ka|^\alpha \xi_\ka$ for some $\alpha>\frac12$.
For $C_b >0$ we then define
$$
\lambda_\ka=|\ka|^2,
\qquad 
B_{\ka,\el,\m}(\chi,\chi')=\frac{\ka\cdot\el^\perp}{|\ka\cdot\el^\perp|}\chi\chi',
\qquad
q_{\ka,\el,\m} = C_b^{-1}\frac{|\ka|^{\alpha-2}
|\ka\cdot\el^\perp|}{|\el|^{\alpha+2}|\m|^\alpha},
$$
for all $\ka$, $\el$, $\m\in\Z^2$ satisfying 
$\ka\cdot\el^{\perp}\neq0$ and $\el+\m=\ka$ (and zero otherwise).
Lemma \ref{alphamajo} ensures that $q_{\ka} < \infty$ and that
$q_\ka \to 0$ as $|\ka| \to \infty$. Taking $C_b$ large enough we have
that $q_\ka < 1$ (note that here we may take $p_\ka=0$). So the recasting 
is complete if we define $ \gamma_\ka=(|\ka|^{\alpha-2}/d_\ka) f_\ka$.
\subsection{A surface growth equation}\label{surface}
The final example illustrates the change in weights needed 
for a higher order equation and the need to consider linear instabilities.
In particular, the linear operator does not generate a diffusion.
Therefore, the Fourier transform is necessary for the stochastic representation.
Consider the following scalar equation arising in some models for surface growth,
$$
\partial_t u=-a_1\Delta^2u-a_2\Delta u-a_3\Delta|\nabla u|^2+a_4|\nabla u|^2+f,
$$
with periodic boundary conditions on $[0,2\pi)^d$, with $d=1,2$ and
$a_i>0$ for $i=1,2,3$. See Raible et al.\ \cite{RaMaLiMoHaSa} for the
derivation of the model, and Bl\"omker et al.\ \cite{DBGuRa} for a
rigorous mathematical treatment using {\pde} techniques. For simplicity,
we assume $a_4=0$ and that the mean value $\int u(t,x)\,dx$ is zero, allowing
us to omit the coefficient $u_\mathbf{0}$.

The equation for the Fourier coefficients is given by
\begin{equation}
\label{e:surf-sys}
\dot u_\ka=-a_1|\ka|^4u_\ka+a_2|\ka|^2u_\ka+a_3|\ka|^2\sum_{\el+\m=\ka}(\el\cdot\m)u_\el u_\m+f_\ka.
\end{equation}
We set $\chi_\ka=|\ka|^\alpha u_\ka$ for $\alpha>0$, with
$\alpha>\max\{d,1 + d/2\}$, and  then choose, for all $\ka \neq 0$,
$$
\lambda_\ka=a_1|\ka|^4,
\qquad
p_\ka = \frac{a_2}{a_1} C_f^{-1} |\ka|^{\alpha-2}
$$
and
$$
B_{\ka,\el,\m}(\chi,\chi')=\frac{\el \cdot \m}{|\el \cdot \m|} \chi \chi',
\qquad
q_{\ka,\el,\m} = C_b^{-1}\frac{a_3|\ka|^{\alpha-2}|\el\cdot\m|}{a_1|\el|^\alpha|\m|^\alpha}
$$
with $B_{\ka,\el,\m}$ and $q_{\ka,\el,\m}$ equal to zero if $\el\cdot\m=0$ or
$\el+\m\neq\ka$. Lemma \ref{alphamajo} guarantees that $p_\ka + q_\ka <1$ when
$C_b,C_f$ are taken large enough and we can define $d_\ka = 1-q_\ka-p_\ka$ and
$ \gamma_\ka=\frac{|\ka|^{\alpha-4}}{a_1d_\ka}f_\ka$ to obtain a system 
in the form \eqref{geneq}.
\subsection{Discussion of extensions and generalisations}\label{ss:moregen}
We now list a
number of possible extensions to the our basic system \eqref{geneq}
for which modified tree representations will hold.
\subsubsection{{\pde}s in general domains with other boundary
conditions} If there is a complete countable set of
$L^2$-eigenfunctions $(e_k)_{k\in\N}$ of $A$ in which to expand
solutions as $u=\sum_{k=1}^\infty u_ke_k$, one can recast
{\pde}s in general domains with various boundary conditions
into a suitable {\ode} setting. This is similar to spectral
Galerkin methods. Consider for instance
$$
\partial_t u=Au+B(u,u)
$$
for a bilinear operator $B$, and suppose that $Ae_k=\lambda_k e_k$.
Then 
$$
\partial_t u_k
=\lambda_ku_k+\sum_{m,l=1}^\infty\langle B(e_m,e_l),e_k\rangle_{L^2}\, u_m u_l.
$$
This can easily be transformed into the general system \eqref{geneq}
by choosing appropriate weights. This would cover our earlier examples, Burgers 
equation, Navier-Stokes or the surface growth equation, in a regular domain with,
for instance, Dirichlet or Neumann boundary conditions. 
Note that \eqref{geneq} is now an $\R$-valued system posed in
$\ell^\infty(\R)$ which is indexed over $\N$ instead of $\Z$.
\subsubsection{Polynomial non-linearities} More general polynomial
non linearities, or several non-linearities, lead to branching
systems where particles split into a larger number of descendants.
Even analytic non-linearities can be handled, with the absolute
values of the power series coefficients controlling the branching
probabilities. 
Note also that a general first order term of the form
$\sum_{\el} p_{\ka,\el} B_{\ka,\el}(\chi_\ka)$, for linear
$B_{\ka,\el}:\C^r \to \C^r$, can be thought of as a branching event
with a single offspring.
This kind of term arises, for example, when the original \pde\ 
contains a multiplication operator $u\mapsto fu$ for a fixed function
$f$.
\subsubsection{Multiplicative forcing} A non-linear forcing term
$F(u,f)$, again with polynomial $F$, can also be recast into a
branching system of {\ode}s. This leads, say in the quadratic case, to
time dependent bilinear operators $B_{\ka,\el,\m}$ whose values depend
on the forcing $\ga_{\ka}(t)$, i.e.\ we obtain terms like a
sum over $q_{\ka,\el,\m} B_{\ka,\el,\m}(\gamma_\el,\chi_{\m})$
in equation \eqref{geneq}.
%%
%%
%%
%%%%%%%%%%%%%%%%%%%%%%%%%%%%%%%%%%%%%%%%%%%%%%%%%%%%%%%%%%%%%%%%%%%%%%%%%%%%%%%%%
%%
%%
\section{The branching particle representation formula} \label{s:rep}
\subsection{Existence and uniqueness}
The next theorem shows that there is a unique local solution to 
\eqref{geneq} taking values in the space $\ell^\infty(\C^r)$ of
bounded families $(a_\ka)_{\ka\in\Z^d}$ of elements of $\C^r$,
with the norm $\|a\|_\infty=\sup_{\ka\in\Z^d}|a_\ka|$,
with $|a_\ka|= \sqrt{a_{\ka} \cdot a_{\ka}^*}$ the norm in $\C^r$. 
We give a simple deterministic proof, but there is also a more
probabilistic proof available (see Corollary \ref{compexun}), in the
spirit of Le Jan and Sznitman \cite{LJS}. 
\begin{theorem}[Unique local existence] \label{locex}
Assume that
$$
\chi(0)\in\ell^\infty(\C^r),\qquad
\ga\in L^{\infty}([0,T],\ell^{\infty}(\C^r))
\quad\text{for all }T>0.
$$                                                                                
Then there exists a time $T_0>0$, depending only on $\chi(0)$, $\gamma$,
and the constants appearing in the equation, such that the mild
formulation \eqref{mildform} has a unique solution
$\chi\in L^{\infty}_\loc([0,T_0),(\C^r)^{\Z^d})$.

Moreover, we have either $T_0=\infty$ or $\|\chi(t)\|_\infty\to\infty$ as $t\to
T_0$. Finally, if the functions $t \mapsto \gamma_\ka(t)$ are $C^k$,
then $t \mapsto \chi_\ka(t)$ are $C^{k+1}$ in time and solve 
equation \eqref{geneq}.
\end{theorem}
\begin{proof}
The proof is a rather standard application of the Banach fixed point
theorem. Let $B$ be a ball of radius $R>0$ centred at the constant function 
with value $\chi(0)$, in the space $L^{\infty}([0,t_\star],\ell^\infty(\C^r))$.
For $\chi \in B$ define $F(\chi)$ by the right-hand side of \eqref{mildform}. 
Then for $R_0=R+\|\chi(0)\|_\infty$,
$$
|F(\chi)_{\ka}(t)-\chi_\ka(0)|
\le \left(\|\chi(0)\|_\infty 
  +C_f p_\ka R_0
  + C_b q_{\ka} R_0^2
  + \|\gamma(t)\|_\infty\right)\!\!(1-\e^{-\lambda_\ka t_*}).
$$
If we choose $R> \|\chi(0)\|_\infty+\sup\|\gamma(t)\|_\infty$
and $t_*$ small enough we see that $F$ maps $B$ into itself. 
Here we have used assumption \eqref{asslocex} to
control the large $|\ka|$s.
Moreover, if $\chi^1$ and $\chi^2$ are in $B$, then for $t \leq t_*$,
$$
|[F(\chi^1)-F(\chi^2)]_{\ka}(t)|
\le(C_f p_\ka +2 C_b R_0 q_{\ka}) 
(1-\e^{-\lambda_\ka t_*}\!\!)\sup_{t \leq t_*}\|\chi^1(t)-\chi^2(t)\|_{\infty}.
$$
Hence $F$ is a strict contraction in $L^{\infty}([0,t_\star],\ell^\infty(\C^r))$
if we choose $t_*$ small enough. Here we need again, for large
$|\ka|$, the assumption \eqref{asslocex}.
                                                                                
The assertion for the time $T_0$ follows in a standard manner
by gluing together local solutions.
The continuity of $t \to \chi_\ka (t) $ follows from the mild
form \eqref{mildform}. It is even differentiable with bounded 
derivative.
The $C^k$-regularity follows by
differentiating \eqref{mildform} and the higher regularity
follows from differentiating \eqref{geneq}.
\end{proof}
A simple global existence result can be proved under the assumptions of 
linear stability and small data.
\begin{proposition}[Global existence for small data] \label{prop:globex}
Under the assumptions of Theorem \ref{locex},
assume that there exists $\delta > 0$ so that
$$
d_\ka |\gamma_\ka| < \delta (1-C_fp_\ka) - C_b \delta^2 q_\ka, 
\quad \mbox{for all $\ka\in\Z^d$.}
$$
Then for each initial condition $\|\chi(0)\|_\infty \leq \delta$,
there is a global solution $\chi$ to equations \eqref{geneq} satisfying 
$ \sup_{t\ge0} \|\chi(t)\|_\infty \leq \delta $.
\end{proposition}
\begin{proof}
Let $v_\ka (t)=|\chi_\ka(t)|$. When $v_\ka(t) \neq 0 $ and $v_\la(t) \leq \delta$ 
for all $\la\in\Z^d$ one 
has the estimate
\begin{equation} \label{temp111}
\partial_t v_\ka \le - \lambda_\ka \left(1-C_fp_\ka\right) v_\ka + 
\lambda_\ka \left(C_b q_{\ka} \delta^2 + d_\ka|\gamma_\ka| \right).
\end{equation} 
The assumption implies that the right hand side of \eqref{temp111} is negative when 
$v_{\ka} = \delta$ and global existence follows from 
a comparison argument for one-dimensional {\ode}s.
\end{proof}
\subsection{The branching tree}
We now give a construction of the branching process that will be used
to represent the solutions of \eqref{geneq}.  We will label particles
of the process with labels taken from the set $\I= \bigcup_{n=0}^\infty
\{0,1,2\}^n$.  The history of a particle $\alpha= (\alpha_1, \ldots,
\alpha_n)$ can be read off  by interpreting $\alpha_j=0$ as the flip
at generation $j$, and $\alpha_j=1$ (or $2$) as being child $1$ (or
$2$) in a binary branching event at generation $j$.

For $\alpha \in \{0,1,2\}^n$ we write $|\alpha|=n$ which we call the
length of the label.  We write $\alpha = \emptyset$ for the single
label of length zero.  When $\alpha=(\alpha_1, \ldots, \alpha_n)$ we
write $\alpha|_{j}$ for the label $\alpha|_j=(\alpha_1, \ldots,
\alpha_j)$ of its ancestor at generation $j \in \{0,1,\ldots,n-1\}$
(and set $\alpha|0= \emptyset$). For $ i \in \{0,1,2\}$ we write
$(i,\alpha)$ for the label  $(i,\alpha_1, \ldots, \alpha_n)$ and
$(\alpha,i)$ for the label  $(\alpha_1, \ldots, \alpha_n,i)$ (or
$(i,\alpha)=(\alpha,i)=(i)$ if $\alpha=\emptyset$).
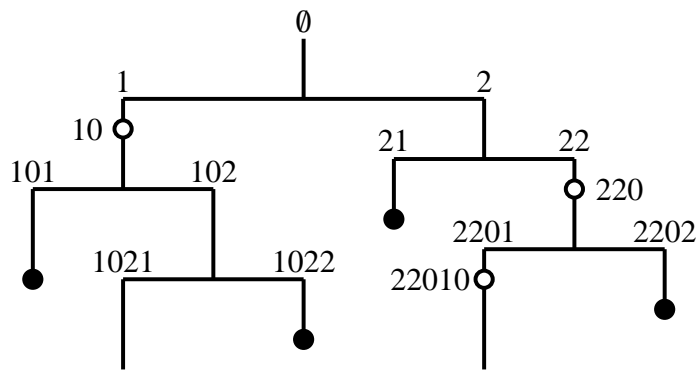
\begin{figure}[h]
\centering
\psset{xunit=1.2cm, yunit=0.4cm, linewidth=1.5pt}
\begin{pspicture}(0,0)(7,12)
\qline(3,11)(3,9)\qline(1,9)(5,9)
\qline(1,9)(1,6)\qline(0,6)(2,6)
\qline(0,6)(0,3)\qline(2,6)(2,3)
\qline(1,3)(3,3)\qline(1,3)(1,0)
\qline(3,3)(3,1)\qline(5,9)(5,7)
\qline(4,7)(6,7)\qline(4,7)(4,5)
\qline(6,7)(6,4)\qline(5,4)(7,4)
\qline(5,4)(5,0)\qline(7,4)(7,2)
\pscircle*(0,3){4pt}\pscircle*(3,1){4pt}
\pscircle*(4,5){4pt}\pscircle*(7,2){4pt}
\pscircle[fillstyle=solid](1,8){4pt}
\pscircle[fillstyle=solid](6,6){4pt}
\pscircle[fillstyle=solid](5,3){4pt}
\rput(3,11.6){$\emptyset$}
\rput(1,9.6){1}\rput(5,9.6){2}
\rput(0.60,8){10}
\rput(0,6.6){101}\rput(2,6.6){102}
\rput(1,3.6){1021}\rput(3,3.6){1022}
\rput(4,7.65){21}\rput(6,7.65){22}
\rput(6.5,6){220}
\rput(5,4.6){2201}\rput(7,4.6){2202}
\rput(4.4,3){22010}
\end{pspicture}
\caption{A tree with branches, deaths ($\bullet$) and flips ($\mathbf{\circ}$).}
\end{figure}
We construct the branching particle systems on a probability space 
equipped with the following independent families of I.I.D. variables: 
$(E_{\alpha})_{\alpha\in\I}$ exponential mean one variables (that
will control the overall rates of branching and flipping);
$(U_{\alpha})_{\alpha\in\I}$ uniform $[0,1]$ variables (that will
control whether a particle flips, branches or dies) and
$((Y_{\alpha}^{(1)}(\ka),Y_{\alpha}^{(2)}(\ka))_{\alpha\in\I,\,\ka\in\Z^d}$
random variables with distribution
$P[Y_{\alpha}^{(1)}(\ka)=\el,\,Y_{\alpha}^{(2)}(\ka)=\m] = q_{\ka,\el,\m}$
(which will control the positions of the two offspring of a particle
that branches).

We now define a system 
$(\hat{\K}_{\alpha}, \tau^B_{\alpha}, \tau^D_{\alpha})_{\alpha\in\I}$
of particle positions, birth and death times,
inductively over the length $n= | \alpha|$ of the labels. 
Fix $\ka \in \Z^d$ and set $\hat{\K}_{\emptyset}=\ka$, $\tau^B_{\emptyset}=0$
and $\tau^D_{\emptyset} = \lambda_{\ka}^{-1} E_{\emptyset}$. Assume that
the positions, birth and death times have been defined for $|\alpha|\leq n$. 
Then, for $\alpha$ of length $n+1$, define birth and death times 
$$
\tau^B_{\alpha} = \tau^D_{\alpha |n},
\qquad
\tau^D_{\alpha}  =  \tau^B_{\alpha} + 
\lambda^{-1}_{\hat{\K}_{\alpha |n}} E_{\alpha},
$$
and the particle positions
$$
\hat{\K}_{\alpha}=\begin{cases}
\hat{\K}_{\alpha|n},                   &\quad\alpha_{n+1}=0,\\
Y^{(1)}_{\alpha}(\hat{\K}_{\alpha|n}), &\quad\alpha_{n+1}=1,\\
Y^{(2)}_{\alpha}(\hat{\K}_{\alpha|n}), &\quad\alpha_{n+1}=2.
\end{cases}
$$
This defines a complete tree of all possible branching and flipping particles 
rooted at $\ka$. In the desired evolution the particles will choose
whether to flip, branch or die according to the probabilities $p_\ka$,
$q_\ka$, $d_\ka$.
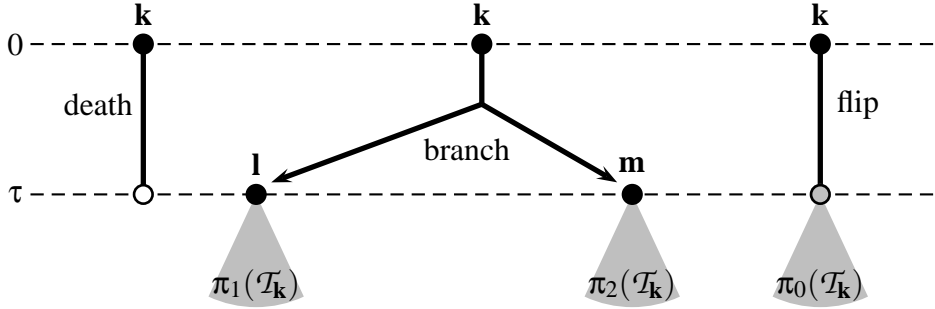
\begin{figure}[h]
\centering
\psset{xunit=1cm, yunit=0.4cm}
\begin{pspicture}(0,0)(12,9)
\psline[linestyle=dashed]{-}(0,3)(12,3)
\psline[linestyle=dashed]{-}(0,8)(12,8)
\rput(-0.2,8){0}
\rput(-0.2,3){$\tau$}
\pswedge*[linecolor=lightgray](3,3){1.5}{245}{295}\rput(3,0){$\pi_1(\T_\ka)$}
\pswedge*[linecolor=lightgray](8,3){1.5}{245}{295}\rput(8,0){$\pi_2(\T_\ka)$}
\pswedge*[linecolor=lightgray](10.5,3){1.5}{245}{295}\rput(10.5,0){$\pi_0(\T_\ka)$}
\psset{linewidth=2pt}
\pscircle*(1.5,8){4pt}\qline(1.5,8)(1.5,3)
\pscircle[fillstyle=solid,linewidth=1pt](1.5,3){4pt}
\rput(0.9,6){death}\rput(1.5,9){$\ka$}
\pscircle*(10.5,8){4pt}\qline(10.5,8)(10.5,3)
\pscircle*[linecolor=lightgray,linewidth=1pt](10.5,3){4pt}
\pscircle[linewidth=1pt](10.5,3){4pt}
\rput(11,6){flip}\rput(10.5,9){$\ka$}
\pscircle*(6,8){4pt}
\rput(5.8,4.5){branch}\rput(6,9){$\ka$}
\qline(6,8)(6,6)
\pscircle*(3,3){4pt}\psline{->}(6,6)(3.2,3.4)\rput(3,4){$\el$}
\pscircle*(8,3){4pt}\psline{->}(6,6)(7.8,3.4)\rput(8,4){$\m$}
\end{pspicture}
\caption{The construction of the tree. At each event time $\tau$, there is
a random selection between either \emph{death}, \emph{flip} or \emph{branch}
of two new particles to states $\el$ and $\m$ (depending on the state $\ka$
of the parent particle)}
\end{figure}

We now define indicator variables $(I_{\alpha})_{\alpha\in\I}$ to
decide whether a particular branch has survived. Define $I_{\emptyset}=1$
and, for $\alpha$ of length $n+1$,
$$
I_{\alpha}=\begin{cases}
1 \quad&\text{if }\alpha_{n+1}=0\text{ and }U_{\alpha}\in[0,p_{\hat{\K}_{\alpha |n}}],\\
1 &\text{if }\alpha_{n+1}\in\{1,2\}\text{ and }U_{\alpha}\in[1- q_{\hat{\K}_{\alpha |n}},1],\\
0 &\text{otherwise.}
\end{cases}
$$ 
Fix an isolated cemetery state $\Delta$ and define, for $|\alpha|=n$,
$$
\K_{\alpha}=\begin{cases}
\hat{\K}_{\alpha} &\text{if }\prod_{j=1}^n I_{\alpha|j} =1,\\
\Delta           &\text{otherwise.}
\end{cases}
$$
The collection 
$\T_{\ka}= (\K_{\alpha},\tau^B_{\alpha},\tau^D_{\alpha})_{\alpha\in\I}$ 
now defines our branching tree rooted at $\ka$. It lives in the space
defined by 
$$
\sT=\left((\Z^d\cup\{\Delta\})\times[0,\infty)\times[0,\infty)\right)^{\I}.
$$
We denote the law of $\T_{\ka}$ on $\sT$ by $\Pb_{\ka}$.

The descendants of any one particle in the tree form a new tree. To 
make this precise we define shift maps $\pi_i:\sT\to\sT$, for $i=0,1,2$ as follows. 
For $\T=(\ka_{\alpha},s_{\alpha},t_{\alpha})_{\alpha\in\I}\in\sT$ we define
a new tree $\pi_i(\T)$ by 
$$
\pi_i(\T) =(\ka_{(i,\alpha)},s_{(i,\alpha)} - t_{\emptyset}, 
t_{(i,\alpha)}- t_{\emptyset})_{\alpha\in\I}.
$$
The tree $\pi_i(\T)$ is meant to be the tree of descendants of the 
particle labelled $(i)$, with their birth and death times shifted so that
particle $(i)$ is born at time $t=0$. 
The construction of the branching particle system from I.I.D. families
implies the following lemma.
\begin{lemma} \label{onestep}
Let $\T_{\ka}=(\K_{\alpha},\tau^B_{\alpha},\tau^D_{\alpha})_{\alpha\in\I}$
have law $\Pb_{\ka}$. Then 
\begin{enumerate}
\item conditional on $\{ \tau^D_{\emptyset} \in ds, \, \K_{(0)} = \ka \}$ the 
tree $\pi_0(\T_{\ka})$ has the law $\Pb_{\ka}$; 
\item conditional on $\{\tau^D_{\emptyset} \in ds, \, \K_{(1)} = \m, \K_{(2)}=\el\}$ 
the trees $\pi_1(\T_{\ka})$ and $\pi_2(\T_{\ka})$ are independent 
and have laws $\Pb_{\m}$ and $\Pb_{\el}$.
\end{enumerate}
\end{lemma}
We want to ensure that the tree has only finitely many branches before time $t$.
Define $N_{[0,t]}:\sT\to\N$ by 
$N_{[0,t]}(\T) = | \{\alpha \in \I: s_{\alpha} \leq t\} |$, that is
the cardinality of the set of particles born before time $t$.
\begin{lemma}\label{alive}
Under $\Pb_{\ka}$ the variables $N_{[0,t]}$ are almost surely finite
for all $t \geq 0$.
\end{lemma}
\begin{proof}
Let $P_\ka(t)= \Pb_{\ka}[N_{[0,t]} < \infty]$. By conditioning on the
values of $\tau^D_{\emptyset}$, $\K_{(0)}$, $\K_{(1)}$, $\K_{(2)}$ and
using Lemma \ref{onestep},
$$
P_\ka(t) = \e^{-\lambda_\ka t} + \int_0^t \lambda_\ka \e^{-\lambda_\ka(t-s)}
  \big[ p_\ka P_\ka(s) + \sum_{\el,\m\in\Z^d} q_{\ka,\el,\m} P_\el(s) P_\m(s) + d_\ka \big] \,ds.
$$
Hence, $(P_\ka(t): \ka\in\Z^d, \, t \geq 0)$ is a bounded, real-valued solution to the equation
\eqref{geneq} with forcing $\ga \equiv 1$ and bilinear operators
$B_{\ka,\el,\m}(\chi, \chi') = \chi \chi'$. By Theorem \ref{locex},
there is only one solution, namely $P_\ka (t) = 1$ for all $\ka, \,t$.
\end{proof}
A simple criterion that ensures that the branching process becomes extinct 
with probability one, that is $\K_{\alpha}= \Delta$ for all large $|\alpha|$, 
is that 
\begin{equation} \label{simplecriterion}
q_\ka\leq d_\ka\quad\text{ and }\quad p_\ka<1\qquad\text{ for all }\ka\in\Z^d.
\end{equation}
Indeed the number of particles alive at time $t$ is
an integer valued process whose successive values, under the
condition \eqref{simplecriterion}, form a sub-critical branching process.
Therefore it eventually reaches zero. The number of values 
$\ka \in \Z^d$ taken by particles before this extinction is almost 
surely finite. The conditions that
$\la_\ka>0$ and $p_\ka<1$ ensure that the extinction time for the
branching particle system is almost surely finite.
Note that, as explained in remark \ref{Remark2.1}, we can always 
choose the system \eqref{geneq} in such a way that 
\eqref{simplecriterion} holds.
\subsection{The evaluation along the tree}
We now fix a forcing functions $\gamma$ and an initial condition
$\chi(0)$. We wish to define evaluation maps $R_t:\sT\to\C^r$ for
$t \geq 0$, which will depend on $\gamma$ and $\chi(0)$. These will
satisfy a recursive property that allows them to be calculated
backwards along the tree.

For the sake of simplicity, we introduce the following abbreviations:
given a branching tree $\T=(\ka_\al,s_\al,t_\al)_{\al\in\I}$ and a
particle labelled $\al\in\I$, with $\ka_\al\neq\Delta$, we say that
the particle has a
\begin{description}
\item[death] if $\ka_{(\al,0)}=\ka_{(\al,1)}=\ka_{(\al,2)}=\Delta$,
\item[flip] if $\ka_{(\al,0)}\neq\Delta$ and
            $\ka_{(\al,1)}$, $\ka_{(\al,2)}=\Delta$,
\item[branch] if $\ka_{(\al,0)}=\Delta$ and $\ka_{(\al,1)}$,
            $\ka_{(\al,2)}\neq\Delta$.
\end{description}
Under each probability $\Pb_\ka$, every particle $\alpha$ for which
$\K_\al\neq\Delta$ must do exactly one of the above three possibilities.
\begin{lemma} 
There exist a family of maps $R_t:\sT\to\C^r$, for $t \geq 0$, satisfying, 
when $N_{[0,t]}(\T) < \infty$, the implicit formula 
\begin{equation}\label{evaluation}
R_t(\T)\!=\!\begin{cases}
\chi_{\ka_{\emptyset}}(0)
  &t_\emptyset \geq t,\\
\gamma_{\ka_{\emptyset}} (t-t_\emptyset)
  &t_\emptyset<t,\text{ death at }\emptyset,\\
C_f\,R_{t-t_{\emptyset}}(\pi_0(\T))
  &t_\emptyset<t,\text{ flip at }\emptyset,\\
C_b\,B_{\ka_{\emptyset},\ka_{(1)},\ka_{(2)}}
      (R_{t-t_\emptyset}(\pi_1(\T)),R_{t-t_\emptyset}(\pi_2(\T))) 
  &\!t_\emptyset<t,\,\text{branch at }\emptyset.
\end{cases}
\end{equation}
\end{lemma}
\begin{proof}
Informally, since the tree is finite when $N_{[0,t]}< \infty$ the value of
$R_t(\T)$ can be calculated backwards along the tree, starting at time $s=t$ 
and working back to time $s=0$:
evaluate the initial
condition $\chi(0)$ at any particles that are alive at time $t$, evaluate the 
forcing function $\ga(s)$ at any particle that dies at time $s <t$, and apply the 
bilinear operators at the times of branching events.

For a careful proof one can define a sequence of approximations $R_{n,t}^\app$ 
in the following way: 
$$
R_{1,t}^\app(\T)=\begin{cases}
\chi_{\ka_\emptyset}(0)
   &\quad\text{if }t_\emptyset\ge t,\\
\gamma_{\ka_\emptyset}(t-t_\emptyset)
   &\quad\text{if }t_\emptyset<t,\text{ death at }\emptyset,\\
1  &\quad\text{otherwise,}
\end{cases}
$$
and $R_{n+1,t}^\app(\T)$ is given by
\begin{equation}\label{evaluation2}
\begin{cases}
\chi_{\ka_{\emptyset}}(0)
  & t_\emptyset\ge t,\\
\gamma_{\ka_{\emptyset}} (t-t_\emptyset)
  &t_\emptyset<t,\text{ death at }\emptyset\\
C_f \, R_{n,t-t_\emptyset}^\app(\pi_0(\T))
  &t_{\emptyset}<t,\text{ flip at }\emptyset\\
C_b \, B_{\ka_{\emptyset},\ka_{(1)},\ka_{(2)}}
\left (R_{n,t-t_\emptyset}^\app(\pi_1(\T)),R_{n,t-t_\emptyset}^\app(\pi_2(\T)) \right)
  &t_\emptyset<t,\text{ branch at }\emptyset.
\end{cases}
\end{equation}
If $N_{[0,t]}<\infty$ then only finitely many iterations are needed and
$R_t(\T) = R_{n,t}^\app(\T)$ for all large $n$.
\end{proof}
In some cases the evaluation can be
written more explicitly.
Let $F(t)$ (respectively $B(t)$) be the number of particles that
have flipped (respectively branched) before time $t$. Let $D(t)$
be the set of labels of particles that have died strictly before time $t$. 

Consider the special case where $r=1$ and that all the 
bilinear forms $B_{\ka,\el,\m}$
coincide with the usual product in $\C$. Then 
the evaluation is given, almost surely under $\Pb_\ka$, by
\begin{equation} \label{exlR}
R_t(\T) = C_b^{B(t)} \, C_f^{F(t)} \,
\prod_{\alpha \in D(t)} \gamma_{\ka_{\alpha}} (t-t_{\alpha}) \,
\prod_{\alpha: t \in [s_{\alpha},t_{\alpha})} \chi_{\ka_{\alpha}}(0). 
\end{equation} 
In the general case, we can only verify, under $\Pb_{\ka}$, that
\begin{equation} \label{exlR2}
|R_t(\T)| \leq  C_b^{B(t)}\,C_f^{F(t)} \, 
\prod_{\alpha \in D(t)} |\gamma_{\ka_{\alpha}} (t-t_{\alpha})| \,
\prod_{\alpha: t \in [s_{\alpha},t_{\alpha})} |\chi_{\ka_{\alpha}}(0)|,
\end{equation} 
and that equality holds in \eqref{exlR2} 
if $|B_{\ka,\el,\m}(\chi,\chi')|=|\chi||\chi'|$ for all
$\ka,\el,\m$ and $\chi,\chi'$.
\subsection{The representation formula}
Consider an initial condition $\chi(0)\in \ell^\infty (\C^r)$,
and a forcing $\gamma\in L^\infty([0,T], \ell^\infty(\C^r))$. The
representation formula for solutions of \eqref{geneq}, 
when the expectation exists, is given by
\begin{equation} \label{repform}
\chi_\ka(t)=\E_\ka \left[ R_t \right], \qquad\ka\in\Z^d.
\end{equation}
\begin{theorem} \label{thm:rep}
Suppose that there exists $C=C(\ga, \chi(0), T)< \infty$ so that
$$
\E_\ka\left|R_t\right|\leq C\qquad\text{for all }\ka\in\Z^d\text{ and all }t\in [0,T].
$$
Then $\chi$ defined in \eqref{repform} is the unique
$L^\infty([0,T],\ell^\infty(\C^r))$ solution of problem \eqref{geneq}
for the data $\ga$, $\chi(0)$.
\end{theorem}
\begin{proof}
Note that uniqueness follows from Theorem \ref{locex}.
Fix a $\ka\in\Z^d$. Conditioning on the values of 
$\tau_{\emptyset}, \K_{(0)}, \K_{(1)}, \K_{(2)}$ and using lemma \ref{onestep}
leads immediately to the mild form of the equation \eqref{mildform}.
The uniform (over $\ka$) integrability is necessary to show that the sum over $\el,\m$
converges.
\end{proof}

In the next two sections we discuss how to check the integrability assumption and 
what to do if it fails. We also see, what happens if the solution fails to be in $\ell^\infty$.
%%
%%
%%
%%%%%%%%%%%%%%%%%%%%%%%%%%%%%%%%%%%%%%%%%%%%%%%%%%%%%%%%%%%%%%%%%%%%%%%%%%%%%%%%%%
%%
%%
%%
\section{The comparison equation}\label{compare}
\subsection{The comparison equation}
The comparison equation for system \eqref{geneq} is formed by 
taking the norm of the data $|\chi_\ka(0)|$ and $|\ga_\ka|$ as new data 
for the system 
\begin{equation}\label{e:comp}
\begin{cases}
\dot{\tilde{\chi}}_\ka
   =\lambda_\ka \left[-\tilde{\chi}_\ka
     +C_f p_\ka \tilde{\chi}_\ka
     +C_b \sum_{\el,\m \in\Z^d} q_{\ka,\el,\m} \tilde{\chi}_{\el} \tilde{\chi}_{\m}
     +d_{\ka} |\gamma_\ka| \right],\\
\tilde{\chi}_\ka (0) = | \chi_\ka (0)|\qquad\text{for }\ka \in \Z^d. 
\end{cases}
\end{equation}
We now look for non-negative real solutions $\tilde{\chi}_\ka(t)$.

We also define a modified evaluation operator $\tilde{R}_t$ on $\sT$ by the 
implicit formula \eqref{evaluation} where we use the new data $|\chi_\ka(0)|$
and $|\ga_\ka|$ and the bilinear operators are replaced by
$\tilde{B}_{\ka,\el,\m}(\chi, \chi') = \chi \chi'$, the normal product of 
real numbers. Then $\tilde{R}_t(\T) \geq 0$ and formally we expect that
\begin{equation} \label{repform2}
\tilde{\chi}_\ka(t)=\E_\ka \left[ \tilde{R}_t\right], \qquad\ka\in\Z^d
\end{equation}
should solve the comparison equation. 

The next theorem confirms this and shows that a finite solution to
the comparison equation \eqref{e:comp} is a sufficient, and sometimes
necessary, condition for the tree expectations $\E_\ka [R_t]$ to exist. 
\begin{theorem}\label{thm:comp}
If the expectations in \eqref{repform2} are finite for all $t \in [0,T]$ and
$\ka\in\Z^d$, then they define a mild solution to the comparison equation
\eqref{e:comp} for which $t\to\tilde{\chi}_\ka(t)$ is continuous on $[0,T]$.

Conversely if there exists a finite mild solution of \eqref{e:comp}, that is
$\tilde{\chi}_{\ka}(t) < \infty$ for $t \in [0,T]$ and $\ka \in \Z^d$, then
the expectations in \eqref{repform2} are finite for $t \in [0,T]$,
and they define the smallest positive solution of \eqref{e:comp}.

Finally, the comparison $\E_\ka[|R_t|] \le \E_\ka[\tilde{R}_t]$
holds, with equality whenever
$ |B_{\ka,\el,\m}(\chi,\chi')|=|\chi||\chi'|$
for all $\ka,\el,\m$ and $\chi,\chi'$. 
\end{theorem}
\begin{proof}
For the first claim of the theorem, condition on the values of
$\tau^D_{\emptyset}$, $\K_{(0)}$, $\K_{(1)}$, $\K_{(2)}$ and apply
Lemma \ref{onestep} to see that the expectations
$\tilde{\chi}_\ka(t)=\E_\ka\left[\tilde{R}_t\right]$ satisfy the
mild form of the comparison equation. Moreover the mild form of
the equation shows  that $e^{\la_\ka t}\,\tilde{\chi}_\ka (t)$
is continuous and increasing in $t$.  Note that the convergence
of the series in the mild formulation  is not a problem here,
because due to positivity, we can use monotone convergence.

For the second part of the theorem, let $\tilde\chi$ be a mild
solution of the comparison equation \eqref{e:comp} in $[0,T]$
with data $|\chi_\ka(0)|$ and $|\ga_\ka|$. Define a sequence
of evaluations on the trees as follows: set
$\tilde{R}_{0,t}^\esp(\T,\tilde\chi)=\tilde\chi_{\ka_{\emptyset}}(t)$
and for each $n \geq 0$,
\begin{equation}\label{marprun}
\tilde{R}_{n+1,t}^\esp(\T,\tilde\chi)=
\begin{cases}
|\chi_{\ka_{\emptyset}}(0)|
  &\!\!\!\!\! t_\emptyset \ge t,\\
|\gamma_{\ka_{\emptyset}}(t-t_\emptyset)|
  &\!\!\!\!\! t_{\emptyset}<t,\text{ death at }\emptyset,\\
C_f \, \tilde{R}_{n,t-t_\emptyset}^\esp(\pi_0(\T),\tilde\chi)
  &\!\!\!\!\! t_\emptyset<t,\text{ flip at }\emptyset,\\
C_b \, \tilde{R}_{n,t-t_\emptyset}^\esp(\pi_1(\T),\tilde\chi)
    \, \tilde{R}_{n,t-t_\emptyset}^\esp(\pi_2(\T),\tilde\chi)
  &\!\!\!\!\! t_\emptyset<t,\text{ branch at }\emptyset.
\end{cases}
\end{equation}
(In the language of next section, the evaluation $\tilde{R}_{n,t}^\esp$
correspond to a \emph{pruning} of the tree after $n$ generations and
the expectation $\E_{\ka}\tilde{R}_{n,t}^\esp(\T,\tilde\chi)$ will solve a 
Picard iteration scheme for \eqref{e:comp}). 

Note that, upon dying, flipping or branching,
particles of length $n$  are evaluated 
using the true solution $\tilde\chi$. Inductively one checks,
by conditioning on the first event, that for all $n\ge0$ 
\begin{equation}\label{marbound}
\E_{\ka}\bigl[\tilde{R}_{n,t}^\esp(\T,\tilde\chi)\bigr] = \tilde\chi_\ka(t).
\end{equation}
Since $N_{[0,t]}<\infty$ under $\Pb_{\ka}$ we have that 
$\tilde{R}_{n,t}^\esp(\T)\to\tilde{R}_t(\T)$ almost surely. By Fatou's lemma
and \eqref{marbound} we find that 
$\E_{\ka}[\tilde{R}_t]\le\tilde{\chi}_{\ka}(t)<\infty$.

The third claim of the theorem is immediate from the upper bound
\eqref{exlR2} and the fact that it is an equality under the
conditions given. 
\end{proof}
\begin{remark}
Note that in the above theorem, and its corollary below,
we do not insist the solutions are bounded in $\ell^\infty$. 

The first two parts of the above theorem show that, when there
exists a finite mild solution $\tilde{\chi}$ to \eqref{e:comp},
the function defined by $\E_{\ka}[\tilde{R}_t]$ is the smallest
solution to \eqref{e:comp} lying below $\tilde{\chi}$.
Note in the case of $\ell^\infty$ solutions there is uniqueness of
solutions, as in Theorem \ref{locex}. 

As in Le Jan and Sznitman \cite{LJS}, it is possible, when
there exists a finite mild solution $\tilde{\chi}$ to \eqref{e:comp}, 
to show that $n\to\tilde{R}_{n,t}^\esp(\T,\tilde{\chi})$ is a
non-negative martingale (with respect to a natural filtration
along generations of the tree). Uniform integrability of this
martingale would then imply that
$\tilde{\chi}_{\ka}(t) = \E_{\ka}[\tilde{R}_{t}]$.
\end{remark}

\begin{corollary}\label{compexun}
Under the conditions of either the first or the second part of
theorem \ref{thm:comp} the expectations $\chi_\ka(t)=\E_\ka [R_t]$
are well defined for $t\in[0,T]$ and $\ka\in\Z^d$ and form a
mild solution to \eqref{geneq}. Moreover, such a solution is
unique among all mild solutions $\chi'$ such that
\begin{equation}\label{othermild}
|\chi'_\ka(t)| \le \E_{\ka}[\tilde{R}_t],\qquad \mbox{for all $\ka\in\Z^d,\ t\in[0,T]$.}
\end{equation}
\end{corollary}
\begin{proof}
The expectations $\E_{\ka}[R_{t}]$ are well defined 
by theorem \ref{thm:comp} as $|R_t| \le\tilde{R}_{t}$.
By conditioning on the first event as before they will solve the mild equation. 
Note that in this case the convergence of the sums in the mild equation is guaranteed 
by the finiteness of the comparison equation.

Let $\chi'$ be a mild solution verifying \eqref{othermild} and define a
sequence of evaluations $R_{n,t}^\esp(\T,\chi')$ for $n\in\N$ as in the
proof of previous theorem, that is $R_{0,t}^\esp(\T)=\chi'_{\ka_\emptyset}(t)$
and, for all $n\ge 1$, $R_{n,t}^\esp(\T)$ is defined as in formula \eqref{marprun}
with data $\chi'$ and $\gamma$ and with products $B_{\ka,\ka_{(1)},\ka_{(2)}}$
in the place of usual product. By assumption \eqref{othermild} and an argument
similar to \eqref{exlR2} it follows that
$$
|R_{n,t}^\esp(\T,\chi')|\le\tilde R_{n,t}^\esp(\T,\tilde\chi),
$$ 
where $\tilde{\chi}_\ka(t) = \E_{\ka}[\tilde R_t]$ and 
$R_{n,t}^\esp(\T,\tilde\chi)$ are taken from the proof of Theorem
\ref{thm:comp}. Moreover, as in that proof, we can show inductively
that $\chi'_\ka(t)=\E_\ka[R_{n,t}^\esp(\T,\chi')]$.

We next note that $R_{n,t}^\esp(\T,\chi') = R_t(\T)$ and
$\tilde R_{n,t}^\esp(\T,\tilde\chi) = \tilde{R}_t(\T)$
on the set $\Omega_{n,t}=\{N_{[0,t]}(\T) \leq n \}$.
Thus,
\begin{align*}
\E_\ka [\tilde{R}_{n,t}^\esp(\T,\tilde\chi)\uno{\Omega^c_{n,t}} ]
& = \tilde\chi_\ka(t) - \E_\ka [ \tilde{R}_{n,t}^\esp(\T,\tilde\chi) \uno{\Omega_{n,t}} ] \\
&=\E_\ka [\tilde R_t] - \E_\ka [ \tilde R_t \uno{\Omega_{n,t}} ]
 = \E_\ka [ \tilde R_t \uno{\Omega^c_{n,t}}],
\end{align*}
and therefore
\begin{align*}
\left|\chi'_\ka(t)-\E_\ka [R_t]\right|
&\le \E_\ka|R_{n,t}^\esp(\T,\chi')-R_t| 
 =  \E_\ka[|R_{n,t}^\esp(\T,\chi')-R_t|\uno{\Omega^c_{n,t}}],\\
&\le \E_\ka[(\tilde R_{n,t}^\esp(\T,\tilde\chi)+\tilde R_t)\uno{\Omega^c_{n,t}}]
 =  2 \E_\ka[\tilde R_t\uno{\Omega^c_{n,t}}].
\end{align*}
Letting $n\to\infty$ we conclude that $\chi'=\chi$, the solution given
by the probabilistic representation.
\end{proof}
%%
%%
%%%%%%%%%%%%%%%%%%%%%%%%%%%%%%%%%%%%%%%%
\subsection{Examples}
We can remove the weights used to cast the equation into our abstract form
and rewrite the comparison equation as equations for the Fourier
coefficients of a scalar {\pde}.

Consider the Burgers equation example discussed in section \ref{burgers}.
Defining $\tilde{u}_\ka = w_\ka^{-1} \tilde{\chi}_\ka$ we obtain 
a comparison equation of the form 
\begin{align*}
&\dot{\tilde{u}}_\ka=-|\ka|^2 \tilde{u}_\ka+\sum_{\el+\m=\ka}
|\m| \tilde{u}_\el \tilde{u}_\m + |f_\ka|.\\
& \tilde{u}_\ka(0) = |u_\ka^0|
\end{align*}
which in the space coordinates corresponds to the scalar equation
$$
\partial_t \tilde{u} = \Delta \tilde{u} + 
\tilde{u} \, (-\Delta)^{\frac12} \tilde{u} + \tilde{f}
$$
where $\tilde{f}$ has Fourier coefficients $|f_\ka|$.
Note that this scalar comparison equation is independent of the choice of
weights (called \emph{majorizing kernels} in Bhattacharya et al.\ \cite{BCDG}).

For the two-dimensional Navier Stokes equation discussed in 
section \ref{navierstokes} the comparison equation for 
$\tilde{\xi}_\ka = |\ka|^{-\alpha} \tilde{\chi}_\ka$ takes the form
$$
\dot{\tilde{\xi}}_\ka = -|\ka|^2 \tilde{\xi}_\ka + \sum_{\el+\m=\ka}
\frac{|\ka \cdot \el^\perp|}{|\el|^2} \tilde{\xi}_\el \tilde{\xi}_\m+|f_\ka|,
$$
which does not have a nice expression in the space variables.

For the surface equations discussed in section \ref{surface} the 
comparison equation becomes 
$$
\partial_t \tilde{u} = 
-a_1 \Delta^2 \tilde{u} -a_2 \Delta \tilde{u} - \Delta|(-\Delta)^{\frac12} \tilde{u}|^2 + \tilde{f},
$$
where the forcing $\tilde{f}$ has Fourier coefficients $|f_\ka|$.

Whenever there is a solution to these scalar comparison equations
with finite Fourier coefficients we obtain
the existence of mild solutions to the corresponding abstract {\ode}s
given by the stochastic representation \eqref{repform} This
in turn is equivalent to the existence of solutions to the original {\pde}s
with finite Fourier coefficients.

However all three scalar comparison equations have quadratic growth and it is
possible to show, 
for example for zero forcing and large enough initial data, that the solutions 
explode in finite time. See for example \cite{LMW} and the references therein
for the case of branching with diffusion.
\begin{remark}
In the case of the $2d$ Navier Stokes,  the $1d$ Burgers, or the
surface equation, the equality in the last part of Theorem \ref{thm:comp} holds. 
This implies that the stochastic representation $\E_\ka[R_t]$ is 
well defined as the expectation of an integrable variable,
if and only if the corresponding 
comparison equation has a solution with finite Fourier coefficients.
In particular, for any suitable weight, the representation will fail to exist 
at the same time, once a Fourier mode in the comparison equation becomes infinite
for all solutions.
\end{remark}
%%
%%
%%
%%%%%%%%%%%%%%%%%%%%%%%%%%%%%%%%%%%%%%%%%%%%%%%%%%%%%%%%%%%%%%%%%%%%%%%%%%%%%%%%%
%%
%%
%%
\section{The pruned approximation} \label{sec4}
\subsection{An {\ode} example of the approximation scheme}
We first explain the main ideas of the approximation scheme on a simple
example, namely the equation $\dot u=-u+u^2$. The solution can be given
by the stochastic representation $u(t)=\E[u(0)^{N_t}]$, where $N_t$ is the
number of particles at time $t$ of a simple rate one branching process
starting from a single particle at time $0$. It's easy to verify that
the representation is well defined for all time $t\ge0$ if and only
if $|u(0)|\le1$, in that the variable $|u(0)|^{N_t}$ becomes
non-integrable for large $t$ when $|u(0)|>1$, while the solutions
of the equation blow up only if $u(0)>1$.

We now give a modification of the branching process. Give each particle 
a label from the integers $\N$. Particles still branch at rate $1$ 
but a particle with label $n$ produces two offspring, one with label $n$ 
and one with label $n-1$. When a particle of type $0$ tries to branch
it simply dies. Start with a single particle with label $n$ and let 
$N_t(n)$ denote the number of particles at time $t$. 
Set $u_n(t)=\E [u(0)^{N_t(n)}]$ for $n \geq 0$ and 
$u_{-1} \equiv 0$. Then $u_n(t)$ solves the following
semi-implicit iterative scheme
$$
\dot u_n = - u_n + u_{n-1}u_n,  \qquad u_n(0)=u(0),\quad\text{for }n \geq 0.
$$
It is straightforward to check that $u_n(t)$ is well defined for all $n$
and $t$. Moreover, $u_n$ converges to the solution $u(t)$ of the original
problem for each initial condition $u(0)\le1$. This yields the stochastic
representation
$$ 
u(t) = \lim_{n \to \infty} \E \bigl[ u(0)^{N_t(n)} \bigr]
$$
valid for all $u(0) \leq 1$ and all $t\ge0$.
\begin{remark}
The seemingly simpler modification
(used by Le Jan and Sznitman \cite{LJS}
for their uniqueness proof and by Bhattacharya et al.\ \cite{BCDG}) 
where a particle with
label $n$ produces two offspring each with label $n-1$, leads to the explicit 
iterative scheme $ \dot u_n = - u_n + u^2_{n-1} $.
Unfortunately, the limit of $u_n(t)$ for large $t$, as $n \to \infty$, 
fails to exist for $u(0) < -1$.  
\end{remark}
The semi-implicit approximation scheme works for other polynomial non-linearities. 
For example, if one considers $\dot u=-u-u^3$, the approximation scheme
$u_n=-u_n-u_{n-1}^2u_n$,
where each particle with label $n$ branches into three particles,
one with label $n$ and two with label $n-1$, is convergent to
the true global solution for any initial condition.
\subsection{A general approximation scheme}
The aim is to define a sequence of approximations $\chi^\pn_\ka(t) $
to our abstract system of {\ode}s \eqref{geneq}. These approximations will 
have a stochastic representation without any integrability problems.

Rather than construct a particle system with labelled particles 
as described in the previous section, we put the modification into the 
evaluation operators. We claim there exists a sequence of evaluation
operators $R_{n,t}:\sT\to\C^r$ satisfying the following implicit relations
on $N_{[0,t]}< \infty$: 
$$
R_{0,t}(\T)=\begin{cases}
\chi_{\ka_{\emptyset}}(0) &\qquad\text{if }t_{\emptyset}\geq t,\\
0 &\qquad\text{otherwise}
\end{cases}
$$
and, for $n\ge1$, $R_{n,t}(\T)$ equals
\begin{equation}\label{evaluation3}
\begin{cases}
\chi_{\ka_{\emptyset}}(0)
  &t_\emptyset\ge t,\\
\gamma_{\ka_{\emptyset}}(t-t_\emptyset)
  &t_\emptyset<t,\text{ death at }\emptyset,\\
C_f \, R_{n,t-t_\emptyset}(\pi_0(\T))
  &t_\emptyset<t,\text{ flip at }\emptyset,\\
C_b \, B_{\ka_{\emptyset},\ka_{(1)},\ka_{(2)}} 
\left(R_{n,t-t_\emptyset}(\pi_1(\T)),
      R_{n-1,t-t_\emptyset}(\pi_2(\T))\right)
  &t_\emptyset<t\text{ branch at }\emptyset.
\end{cases}
\end{equation}
The existence of $R_{n,t}$ can be established exactly 
as in Lemma \ref{evaluation}. The intuitive link with the 
labelled particle picture in the last section is that
$R_{n,t}(\T)$ corresponds to the evaluation operator
applied to the tree started at a particle with label $n$
at position $\ka$. 
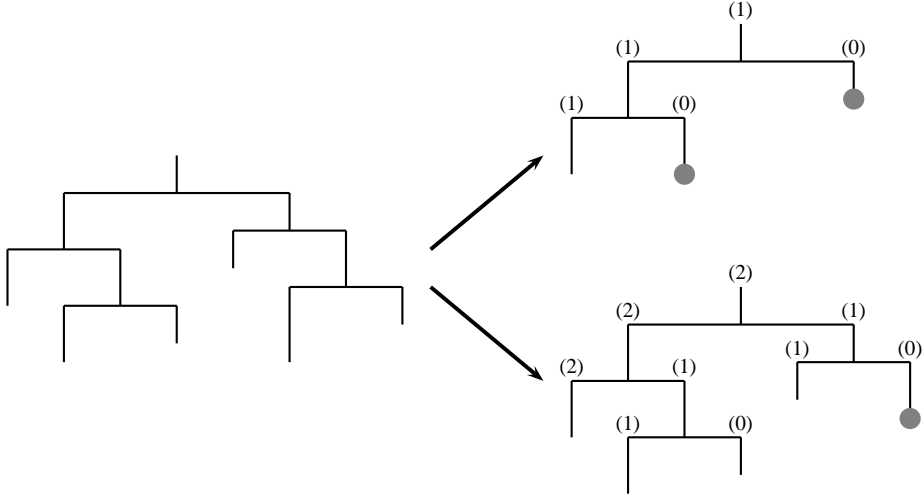
\begin{figure}[h]\label{fig:prune}
\centering
\psset{xunit=0.75cm, yunit=0.25cm}
\begin{pspicture}(0,-7)(16,18)
\qline(3,11)(3,9)\qline(1,9)(5,9)
\qline(1,9)(1,6)\qline(0,6)(2,6)
\qline(0,6)(0,3)\qline(2,6)(2,3)
\qline(1,3)(3,3)\qline(1,3)(1,0)
\qline(3,3)(3,1)\qline(5,9)(5,7)
\qline(4,7)(6,7)\qline(4,7)(4,5)
\qline(6,7)(6,4)\qline(5,4)(7,4)
\qline(5,4)(5,0)\qline(7,4)(7,2)
\psline[linewidth=1.5pt]{->}(7.5,4)(9.5,-1)
\psline[linewidth=1.5pt]{->}(7.5,6)(9.5,11)
\qline(13,4)(13,2)\qline(11,2)(15,2)
\rput(13,4.7){\tiny(2)}\rput(11,2.7){\tiny(2)}\rput(15,2.7){\tiny(1)}
\qline(11,2)(11,-1)\qline(10,-1)(12,-1)
\rput(10,-0.3){\tiny(2)}\rput(12,-0.3){\tiny(1)}
\qline(10,-1)(10,-4)\qline(12,-1)(12,-4)
\qline(11,-4)(13,-4)\qline(11,-4)(11,-7)
\rput(11,-3.3){\tiny(1)}\rput(13,-3.3){\tiny(0)}
\qline(13,-4)(13,-6)\qline(15,2)(15,0)
\qline(14,0)(16,0)\qline(14,0)(14,-2)
\rput(14,0.7){\tiny(1)}\rput(16,0.7){\tiny(0)}
\qline(16,0)(16,-3)
\pscircle*[linecolor=gray](16,-3){4pt}
\qline(13,18)(13,16)\qline(11,16)(15,16)
\rput(13,18.7){\tiny(1)}\rput(11,16.7){\tiny(1)}\rput(15,16.7){\tiny(0)}
\qline(11,16)(11,13)\qline(10,13)(12,13)
\rput(10,13.7){\tiny(1)}\rput(12,13.7){\tiny(0)}
\qline(10,13)(10,10)\qline(12,13)(12,10)
\qline(15,16)(15,14)
\pscircle*[linecolor=gray](15,14){4pt}
\pscircle*[linecolor=gray](12,10){4pt}
\end{pspicture}
\caption{The tree on the left can be pruned by starting, top right,
with a particle labelled $(1)$, or, bottom right, by a particle
labelled $(2)$. Circles mark the death via pruning. If a label
larger than $(2)$ is given to the starting particle, the tree
is unpruned.}
\end{figure}

The implicit relation implies that if 
$N_{[0,t]} \leq m$ then $R_{n,t}(\T) = R_t(\T)$ whenever
$n \geq m$. Moreover when $R_{n,t}(\T) \neq R_t(\T)$ then
$R_{n,t}(\T)=0$. Thus there exist increasing sets
$\Omega_{n,t}\subset\sT$ so that
\begin{equation} \label{singular}
R_{n,t}(\T_\ka) = R_t(\T_\ka) \, \uno{\Omega_{n,t}}
\quad\text{and}\quad
\{ N_{[0,t]} <\infty\}\subseteq\bigcup_n\Omega_{n,t}.
\end{equation}

We now define the stochastic representation using these modified
evaluations by 
\begin{equation} \label{repapprox}
\chi^\pn_\ka(t) = \E_\ka \left[ R_{n,t} \right].
\end{equation}
The fact that this expectation is always well defined is part of
the following result.
\begin{proposition}\label{prop:approx}
Suppose that $\chi(0)\in\ell^{\infty}(\C^r)$ and
$\ga \in L^{\infty}([0,T],\ell^{\infty}(\C^r))$.
Then the expectations in \eqref{repapprox} are well defined and 
$\chi^\pn_\ka(t)$  are the unique $L^{\infty}([0,T],\ell^{\infty}(\C^r))$ 
mild solution to the following approximation scheme
\begin{align}\label{approxscheme}
&\dot\chi^\pz_\ka = - \lambda_\ka \chi^\pz_\ka,\\
&\dot\chi_\ka^\pn = \lambda_\ka \bigl[-\chi_\ka^\pn
  + C_f p_\ka\chi_\ka^\pn
  + C_b \sum_{\el,\m\in\Z^d}q_{\ka,\el,\m}
   B_{\ka,\el,\m}(\chi_{\el}^\pn,\chi_{\m}^\pnm)
  + d_\ka \gamma_\ka\bigr],\notag
\end{align}
with initial condition $\chi^\pn_\ka(0) = \chi_\ka(0)$ 
for all $\ka\in\Z^d$ and $n\in\N$. 
\end{proposition}
\begin{proof}
The local existence and uniqueness of solutions for the approximation scheme, 
follows from the same methods as in the proof of Theorem \ref{locex},
plus an inductive argument in $n \geq 0$. The fact that solutions are globally 
defined follows, again by induction, from the simple estimate
\begin{multline*}
|\chi^\pn_\ka(t)|\le\|\chi(0)\|_\infty + \sup_{t\in[0,T]}\|\gamma\|_\infty+\\
  + \lambda_\ka\int_0^t \e^{-\lambda_\ka(t-s)} 
    \bigl( C_f + C_b \|\chi^\pnm\|_\infty\bigr)\|\chi^\pn\|_\infty\,ds,
\end{multline*}
which, using induction and Gronwall's lemma, easily gives boundedness of 
$\|\chi^\pn\|_\infty$ in each interval $[0,T]$.

In order to prove that the stochastic representation \eqref{repapprox}
is well defined, we use a comparison argument, as in Section \ref{compare}. The
comparison equation for the approximation scheme is given by
$$
\dot{\tilde{\chi}}_\ka^\pn = \lambda_\ka \bigl[
-\tilde{\chi}_\ka^\pn + p_\ka C_f \tilde{\chi}_\ka^\pn + 
C_b \sum_{\el,\m\in\Z^d} q_{\ka,\el,\m} 
\tilde{\chi}_\el^\pn \tilde{\chi}_\m^\pnm + d_\ka |\gamma_\ka|\bigr]
$$
and the evaluation $\E_\ka|R_{n,t}|$ is finite
as long as the $\tilde{\chi}_\ka$ are finite. 
But this follows by the same arguments as in first part of this proof. 
Again $\E_\ka|R_{n,t}|\le\tilde{\chi}_\ka\le C$ 
for all $\ka\in\Z^d$ and all $t\in[0,T]$ with constant $C$ depending
only on $T$, $\chi(0)$, and $\gamma$. 

Finally, the expectations $\E_\ka[R_{n,t}]$ do form the unique solution to
the approximation scheme by conditioning on the first branch of the tree
as in Theorem \ref{thm:rep}.
\end{proof}

In the integrable case, that is where $\E_\ka|R_t|< \infty$, we have 
immediately from \eqref{singular} that
$$
\lim_{n \to \infty} \E_\ka[R_{n,t}] = \E_\ka[R_t].
$$
In particular, when the expectations $ \E_\ka[R_t]$ are bounded over 
$t \in [0,T]$ and $\ka \in \Z^d$ this implies the solutions of the 
approximation scheme converge to those of the original system 
\eqref{geneq}. Our interest, however, is in the non-integrable case and 
we aim to show that convergence of the approximation scheme directly and 
deduce that the limit $\lim_{n \to \infty} \E_\ka[R_{n,t}]$ exists
and defines a stochastic representation for all times $t>0$.
%%
%%
%%
%%%%%%%%%%%%%%%%%%%%%%%%%%%%%%%%%%%%%%%%%%%%%%%%%%%%%%%%%%%%%%%%%%%%%%%%%%%%
%%
\subsection{Global convergence of the stochastic approximation} 
The aim of this section is to give a few details of one
example where the approximation scheme defined by the 
pruned representation converges, even when the
direct stochastic representation fails to be 
integrable. In contrast to the previous section, we 
use {\pde} methods. The convergence depends crucially
on the equation and how the pruning is done, as not
all approximation schemes will converge globally.

For simplicity we work with the one dimensional Burgers
equation with forcing \eqref{burgers}. In Subsection
\ref{suse:burg} we recast the equation 
into our abstract form by considering the 
weighted Fourier coefficients
$$
\chi_\ka(t) = w_\ka u_{\ka}(t)
$$ 
where, as in section \ref{suse:burg}, the weights are given by
$w_\ka=(1\vee|\ka|^\alpha)$ for some $\alpha >1$.
If we assume the Fourier coefficients of the initial condition
satisfy
\begin{equation} \label{hypburg1}
\sup_{\ka} \{|u_\ka(0)|\,w_\ka \} < \infty 
\end{equation}
and the forcing function $f$ satisfies
\begin{equation} \label{hypburg2}
\sup_{\ka} \sup_{t \in [0,T]} \{|f_\ka(t)|\,w_\ka\} < \infty, 
\qquad \text{for all }T>0, 
\end{equation}
then proposition \ref{prop:approx} implies there is a unique global
solution $\chi^\pn_\ka(t)$,
given by \eqref{repapprox}, to the approximation equations
\eqref{approxscheme}.
\begin{theorem}\label{thm:prubur}
Assume, in addition to \eqref{hypburg1} and \eqref{hypburg2}, that 
$u(0)\in H^1$ and  $f\in L^\infty_\loc([0,\infty),L^\infty)$.
Consider the pruned approximation of the previous section. Then the limit 
\begin{equation}\label{e:limB}
\chi_\ka(t)
=\lim_{n\to\infty} \chi_\ka^\pn(t) = \lim_{n \to \infty}
\E_\ka \left[ R_{n,t} \right]
\end{equation}
exists for all $t \geq 0$ and all $\ka\in\Z$ 
and defines a global solution of the Fourier-transformed 
Burgers equation.
\end{theorem}
\begin{proof}
Define
$$
u_\ka^\pn(t) = w_\ka^{-1}\chi_\ka^\pn(t).
$$
Since $\chi^\pn$ is bounded we may reconstruct from these coefficients
the function
$$
u^\pn(t) = \sum_{\ka \in \Z} u_{\ka}^\pn(t)\e^{ \im k x}.
$$
Using the representation of Proposition \ref{prop:approx},
we see that, on the level of {\pde}s, $u^\pn$ solves
the approximation scheme given by 
$$
\begin{cases}
\partial_t u^\pz=\partial_x^2 u^\pz,\\
\partial_t u^\pn=\partial_x^2 u^\pn + \partial_x u^\pn \, u^\pnm+f,\\
u^\pn(0)=u(0).
\end{cases}
$$

Fix $T>0$ and set $C(f,T)=\sup_{t\in[0,T]}\|f(t)\|_{L^\infty}$.
We first use a maximum principle argument to show  
\begin{equation}\label{e:maxprinc}
\sup_{t\in[0,T]}\|u^\pn(t)\|_{L^\infty}\le \|u(0)\|_{L^\infty}+C(f,T) \, T
\qquad\text{for all }n\in\N.
\end{equation}
We now derive an a priori estimate for the solution. 
The following calculation applies to sufficiently smooth 
functions and standard approximation techniques imply that 
the resulting bound holds for the solutions above. 
Using \eqref{e:maxprinc}, we find
\begin{align*}
\frac12\frac{d}{dt}\|\partial_x u^\pn\|^2_{L^2}
& =    -\|\partial_x^2 u^\pn\|^2_{L^2}
       -\int_0^{2\pi}\partial_xu^\pn u^\pnm\partial_x^2u^\pn\,dx\\
&\quad -\int_0^{2\pi}\partial_x^2 u^\pn \,f\,dx\\
&\le -\| \partial_x^2u^\pn\|^2_{L^2}
      + C\| \partial_x^2u^\pn\|^{3/2}_{L^2}
      + C\| \partial_x^2u^\pn\|_{L^2}\;,
\end{align*}
where we have used the Poincar\'e and Cauchy-Schwartz inequalities. 
Note that the constant $C>0$ depends only on $T$, $C(f,T)$,
and $u(0)$. Thus we find another constant, also denoted $C$,  such that 
for all $n\in\N$
$$
\sup_{t\in[0,T]}\|u^\pn(t)\|^2_{H^1}\le C, 
\qquad \int_0^T\|u^\pn(t)\|^2_{H^2}\,dt \le C 
$$
and
$$
\int_0^T\|\partial_t u^\pn(t)\|^2_{L^2}\,dt\le C.
$$
We now use standard methods to show that we have a solution
of the limiting equation (cf. for example Temam \cite{Tem}).
Indeed by  compactness results, there is a subsequence
$(n_k)_{k\in\N}$, such that $u^{n_k}\to u$ weakly
in $L^2([0,T],H^2)$ and $H^1([0,T],L^2)$, and strongly in 
$L^p([0,T],L^2)$ for any $p>1$.
Thus $u$ is the weak solution of Burgers equation,
i.e. it solves the {\pde} in $L^2([0,T],L^2)$.
As weak solutions of the Burgers equation are unique, we can neglect
the subsequence, as any limiting point of $u^\pn$ 
defines the same solution $u$.
Finally, the convergence is strong enough, in order to have all
Fourier coefficients convergent. Thus  for all $\ka\in\Z$
the Fourier coefficients $u_\ka$ of $u$ are given by
$$
u_\ka(t)
=\lim_{n \to\infty} u_\ka^\pn(t)
=\lim_{n \to\infty}w_\ka^{-1} \chi_\ka^\pn(t).
$$
\end{proof}
\begin{remark}
We point out that the assumptions of the 
previous theorem are by no means optimal.
We have used a simplified method of proof, in order to 
provide an example in a simple context.
In particular the constraint on the  initial condition can be relaxed.
Furthermore, using regularisation properties of the {\pde}, we can always get
sufficiently smooth initial conditions, if we wait a small amount of
time. 
\end{remark}
%%
%%
%%
%%%%%%%%%%%%%%%%%%%%%%%%%%%%%%%%%%%%%%%%%%%%%%%%%%%%%%%%%%%%%%%%%%%%%%%%%%%%%%%%%
%%

%%
%%
%%
\end{document}